\documentclass[letterpaper, 10 pt, conference, twocolumn]{ieeeconf}

\IEEEoverridecommandlockouts                              

\usepackage{graphics} 
\usepackage{epsfig} 

\usepackage{times} 
\usepackage{amsmath} 
\usepackage{amssymb}  
\usepackage{dsfont}
\usepackage{booktabs}
\usepackage{mathtools}
\usepackage{subcaption}
\usepackage{soul}

\newtheorem{definition}{Definition}

\newtheorem{prop}{Proposition}

\usepackage{nicefrac}       
\usepackage{microtype}      
\usepackage{tabto}
\usepackage{xcolor}

\usepackage{enumerate}
\usepackage{caption}
\usepackage{subcaption}
\usepackage{graphicx}
\usepackage{stfloats}

\usepackage[notocbasic]{nomencl}
\makenomenclature
\usepackage{etoolbox}

\usepackage{algorithm}
\usepackage[noend]{algorithmic}

\usepackage[colorlinks=true,linkcolor=blue]{hyperref}       
\usepackage{url}            

\usepackage{xpatch}
\xpatchcmd{\thenomenclature}{%
  \section*{\nomname}
}{
}{\typeout{Success}}{\typeout{Failure}}

\newcommand{\norm}[1]{\left\lVert {#1} \right\rVert}

\newcommand{\rlpar}[1]{\left ( {#1} \right)}

\title{\LARGE \bf
A Policy Iteration Algorithm for N-player General-Sum Linear Quadratic Dynamic Games
}

\author{Yuxiang Guan \quad \quad Giulio Salizzoni \quad\quad Maryam Kamgarpour \quad \quad Tyler H. Summers
\thanks{Y. Guan and T. Summers are with the Control, Optimization, and Networks Lab, University of Texas at Dallas. M. Kamgarpour and G. Salizzoni are with the Systems Control and Multi-Agent Optimization Research Lab, EPFL. This work was supported by the United States Air Force Office of Scientific Research under Grant FA9550-23-1-0424 and the National Science Foundation under Grant ECCS-2047040.}
}

\begin{document}

\maketitle
\thispagestyle{empty}
\pagestyle{empty}

\begin{abstract}
We present a policy iteration algorithm for the infinite-horizon N-player general-sum deterministic linear quadratic dynamic games and compare it to policy gradient methods. We demonstrate that the proposed policy iteration algorithm is distinct from the Gauss-Newton policy gradient method in the N-player game setting, in contrast to the single-player setting where under suitable choice of step size they are equivalent. We illustrate in numerical experiments that the convergence rate of the proposed policy iteration algorithm significantly surpasses that of the Gauss-Newton policy gradient method and other policy gradient variations. Furthermore, our numerical results indicate that, compared to policy gradient methods, the convergence performance of the proposed policy iteration algorithm is less sensitive to the initial policy and changes in the number of players.
\end{abstract}

\section{Introduction}
In recent years, the field of multi-agent reinforcement learning (MARL) has attracted significant interest within the reinforcement learning (RL) community. This interest has led to a series of successful developments in approximately solving sequential multi-agent decision-making problems, such as multi-robot control \cite{matignonCoordinated2012}, autonomous driving \cite{shalev-shwartzSafe2016}, and the networking of communication packages \cite{luong_applications_2019}. 
Despite these successes, a comprehensive theoretical understanding of how MARL algorithms perform in environments where cooperation and/or competition among agents exists remains elusive. Recently, there has been a surge of interest in analyzing the performance of policy gradient algorithms within the context of linear quadratic dynamic games (LQDGs). LQDGs present a compelling framework for evaluating the efficacy of MARL algorithms in continuous state and action spaces due to their ability to admit a Nash equilibrium in linear feedback policies. Moreover, this equilibrium can be determined by solving a set of coupled Riccati equations. Properties of this equilibrium have been thoroughly analyzed in \cite{engwerdaScalar1998,basar1998,possieri_algebraic_2015,basarUniqueness1976,lukesGlobal1971}. 

A majority of existing literature \cite{buGlobal2019,zhangPolicy2019,zhangDerivativefree2021} has emphasized zero-sum LQDGs with two players, where it has been shown that certain policy gradient methods have global convergence guarantees in such settings. However, some negative results in \cite{mazumdar2020policy} suggested that the vanilla/standard policy gradient method has no guarantees of even local convergence in infinite-horizon general-sum deterministic LQDGs. This result indicates potential limitations for gradient-type methods in such games. The natural policy gradient and vanilla/standard policy gradient methods demonstrated global convergence in a finite-horizon stochastic LQDG \cite{hambly2023policy}, given appropriate step size and the introduction of specific noise levels into the dynamic system. However, the natural policy gradient method may fail to converge to the Nash equilibrium of a deterministic LQDG without careful selection of initial policies and step sizes.

The policy iteration algorithm is well-known in single-player  settings for computing optimal policies. It comprises two components: policy evaluation and policy improvement. This algorithm has been extensively studied in dynamic programming and RL and bears a close relationship to the Newton method \cite{puterman1979convergence,bertsekas2022lessons}. In the context of single-player LQDGs, the standard policy iteration algorithm \cite{howard:dp} is equivalent to an instance of the Gauss-Newton policy gradient method with a specific step size. 
Several policy iteration-based RL algorithms have been studied for solving multi-player nonzero-sum differential games \cite{vamvoudakis2017game}. Recently, policy iteration algorithms have been developed for solving N-player nonzero-sum LQDGs \cite{nortmannNash2024} which explicitly formulate the policy evaluation and policy update steps. An off-line policy iteration-based RL algorithm was introduced for a two-player nonzero-sum LQDG in \cite{yangDataDriven2020}. This algorithm is a special two-player case of the one we are presenting here. However, a more general policy iteration algorithm for N-player general-sum LQDGs, along with a comparison of its convergence performance against policy gradient methods, has not yet been undertaken.

In this paper, our main contribution is to present a policy iteration algorithm for the infinite-horizon N-player general-sum deterministic LQDGs and compare it to policy gradient methods. 
In contrast to the single-player setting, where the proposed policy iteration algorithm and the Gauss-Newton policy gradient method are equivalent under suitable choice of step size, we show that they are not equivalent in the N-player setting. We illustrate in numerical experiments that the convergence rate of the proposed policy iteration algorithm significantly surpasses that of the Gauss-Newton policy gradient method and other policy gradient variations. Furthermore, our numerical results indicate that, compared to policy gradient methods, the convergence performance of the proposed policy iteration algorithm is less sensitive to the initial policy and changes in the number of players.

The rest of the paper is organized as follows. In Section~\ref{sec:problem_formulation}, we present the formulation of the infinite-horizon N-player general-sum deterministic LQDGs. Section~\ref{sec:policy_iteration} provides a detailed description of the proposed policy iteration algorithm, which is designed to solve these games. The distinction between the proposed policy iteration algorithm and the Gauss-Newton policy gradient method within the N-player game setting is elucidated in Section~\ref{sec:policy_gradient}. Section~\ref{sec:numerical_experiments} showcases the results of our numerical experiments. Finally, we conclude our findings and discussions in Section~\ref{sec:conclusions}.


\section{Problem Formulation: N-player General-Sum Deterministic LQDGs with Infinite-Horizon} \label{sec:problem_formulation}
We consider a discrete-time N-player general-sum deterministic LQDG over an infinite-horizon with dynamics
\begin{equation}\label{eqn:dyna_sys}
    x_{t+1} = A x_t + \sum_{i=1}^N B^i u_t^i,
\end{equation}
where $x_t \in \mathbb{R}^n$ is the system state with the initial value $x_0$ drawn from a Gaussian distribution with $\mathbb E[x_0] = 0$ and $\mathbb E[x_0 x_0^\top]$ = $X_0$, $u_t^i \in \mathbb{R}^{m_i}$ is the control input of player $i = 1,\ldots,N$, and $A \in \mathbb{R}^{n \times n}$ and $B^i \in \mathbb{R}^{n \times m_i}$ are referred to as system matrices. Each player's objective is to minimize their infinite-horizon cost function
\begin{equation}
    \min_{\{u_t^i\}_{t=0}^\infty} \mathbf{E}_{x_0} \left[ \sum_{t=0}^\infty x_t^\top Q^i x_t + (u_t^i)^\top R^i u_t^i \right],
\end{equation}
where $Q^i \in \mathbb{R}^{n \times n}$ and $R^i \in \mathbb{R}^{m_i \times m_i}$ are symmetric matrices that parameterize the quadratic stage costs.


We consider a memoryless perfect state information structure for all the players. That is, each player $i$ seeks a stationary linear feedback policy of the form $u_t^i = K^i x_t$ that minimizes their cost. The policies of all players can be specified by a set of gain matrices $\mathbf{K} = ( K^1, \ldots, K^N )$. Player $i$'s cost induced by the joint policy $\mathbf{K}$ is given by
\begin{equation}
    J^i(\mathbf{K}) \coloneqq \mathbf E_{x_0} \left[ \sum_{t=0}^\infty \left( x_t^\top Q^i x_t + (K^ix_t)^\top R^i (K^ix_t)\right) \right].
\end{equation}
A standard solution concept for general-sum games is a Nash equilibrium. At a Nash equilibrium, no player can unilaterally improve their cost by deviating from their equilibrium policy, defined as follows:
\begin{definition}
    A stationary linear feedback Nash equilibrium for an infinite-horizon general-sum deterministic LQDG is a collection of policies $\mathbf{K}^* = (K^{1*},\ldots,K^{N*})$ such that:
    \begin{equation*}
        J^i(K^{1*},\ldots,K^{i*},\ldots,K^{N*}) \leq J^i(K^{1*},\ldots,K^{i},\ldots,K^{N*}), 
    \end{equation*}
    for each player $i = 1,\ldots,N$.
\end{definition}
Our goal is to study and compare various algorithms for computing Nash equilibria for general-sum LQDGs.



\section{Policy Iteration for N-player General-Sum Deterministic LQDGs with Infinite-Horizon} \label{sec:policy_iteration}
In this section, we first present the well known value iteration algorithm for computing a Nash equilibrium of the infinite-horizon general-sum LQDGs in Section \ref{subsec:value_iteration}. Then we propose a policy iteration algorithm as an alternative to the value iteration algorithm for solving the general-sum LQDGs in Section \ref{subsec:policy_iteration}. 

\subsection{Value Iteration to Compute Nash Equilibrium Policies}\label{subsec:value_iteration}
Value iteration is a standard algorithm for computing feedback Nash equilibrium policies and cost functions in dynamic games. It utilizes principles from dynamic programming for optimal control and is discussed extensively in \cite{basar1998}. Initializing the cost function parameters for each player with $P_0^i = Q^i$, the value iteration algorithm updates the cost parameters and corresponding policies for each player $i = 1,...,N$ at iteration $k = 0,1,...$ via
\begin{equation}\label{eqn:policy_gain}
    K_{k}^{i} = -\left( R^i + \rlpar{B^i}^\top P_{k}^{i} B^i \right)^{-1} \rlpar{B^i}^\top P_{k}^{i} \overline{A}_k^{i},
\end{equation}
\begin{align}\label{eqn:value_fun}
    \begin{split}
        P_{k+1}^{i} = Q^i +& \rlpar{K_{k}^{i}}^\top R^i K_{k}^{i} + \rlpar{\overline{A}_k^{i} + B^i K_{k}^{i}}^\top \\ &P_{k}^{i} \rlpar{\overline{A}_k^{i} +  B^i K_{k}^{i}},
    \end{split}
\end{align}
where $\overline{A}_k^i = A + \sum_{j=1,j \neq i}^N B^j K_k^j$. If these iterations converge, then the limiting policies $K^{i*} = \lim_{k\rightarrow \infty} K_k^i$ and cost parameters $P^{i*} = \lim_{k\rightarrow \infty} P_k^i$ satisfy
\begin{equation}\label{eqn:policy_gain_ss}
    K^{i*} = -\left( R^i + (B^i)^\top P^{i*} B^i \right)^{-1} (B^i)^\top P^{i*} \overline{A}^{i*},
\end{equation}
\begin{align}\label{eqn:value_fun_ss}
    \begin{split}
        P^{i*} = Q^i +& \rlpar{K^{i*}}^\top R^i K^{i*} + \rlpar{\overline{A}^{i*} + B^i K^{i*}}^\top \\ &P^{i*} \rlpar{\overline{A}^{i*} +  B^i K^{i*}},
    \end{split}
\end{align}
a set of coupled algebraic Riccati equations for each player $i = 1,...,N$. Then $K^{i*}$ and $P^{i*}$ are Nash equilibrium policies and cost function parameters for the infinite-horizon general-sum LQDG. The following result provides a sufficient condition for convergence to such a Nash equilibrium.

\begin{prop} [Proposition 6.3 from \cite{basar1998}] \label{prop:LQDG_VI}
    Suppose the above value iteration \eqref{eqn:policy_gain} and \eqref{eqn:value_fun} converges to  $\{ K^{i*}, P^{i*}, \ i \in N \}$ which satisfy \eqref{eqn:policy_gain_ss} and \eqref{eqn:value_fun_ss}, and further suppose that for each $i\in N$ the pair $\big(A + \sum_{j=1,j \neq i}^N B^j K^{j*}$, $B^i\big)$ is stabilizable and the pair $\big( A + \sum_{j=1,j \neq i}^N B^j K^{j*}$, $Q^i + (K^{i*})^\top R^i K^{i*} \big)$ is detectable. Then stationary feedback policies $u_t^{i*} = K^{i*}x_t$ provide a Nash equilibrium solution for the infinite-horizon general-sum LQDGs, leading to the finite infinite-horizon Nash equilibrium cost $x_0^\top P^{i*} x_0$ for player $i$. 
\end{prop}
Under an additional assumption that stage cost parameters satisfy $Q^i \succeq 0$ and $R^i \succ 0$ for each player, the value iteration expressions on the right side of \eqref{eqn:policy_gain} and \eqref{eqn:value_fun} are unique (and correspond to linear feedback Nash equilibrium policies for finite-horizon LQDGs for fixed values of the value iteration index $k$). If the value iteration algorithm converges, the corresponding unique limiting policies $K^{i*} = \lim_{k\rightarrow \infty} K_k^i$ and cost parameters $P^{i*} = \lim_{k\rightarrow \infty} P_k^i$ provide a Nash equilibrium solution for the infinite-horizon general-sum LQDGs. However, the coupled algebraic Riccati equations \eqref{eqn:value_fun_ss} may admit other solutions that are not related to solutions of corresponding finite-horizon dynamic games from the value iteration algorithm \cite{basar1998,papavassilopoulos_linear-quadratic_1984}.



Policy iteration is another dynamic programming-based algorithm that can be utilized to solve general-sum games by iteratively computing a solution to the coupled algebraic Riccati equations \eqref{eqn:value_fun_ss}. We are interested in comparing the convergence properties of policy iteration with variations of policy gradient algorithms for solving general-sum LQDGs. In particular, we aim to study whether policy iteration and policy gradient algorithms converge to the same Nash equilibrium as value iteration or perhaps other Nash equilibria (if they exist), and their respective convergence rates.

\subsection{Proposed Policy Iteration Algorithm}\label{subsec:policy_iteration}
Algorithm \ref{alg:policy_iteration_algo} presents the proposed policy iteration algorithm for the infinite-horizon general-sum LQDGs. Policy iteration has been extensively studied for computing optimal policies for (single-player) optimal control problems, including linear quadratic problems \cite{kleinman1969stability,hewer1971iterative}. Policy iteration begins with an initial stabilizing policy and iterates on two main steps, which are analogous to single-player setting: (1) policy evaluation, which computes the expected costs under the current policy; and (2) policy update, which updates policy under the current cost functions. The policy evaluation step corresponds to solving a set of coupled Lyapunov equations \eqref{policyevaluation}, which relates to the Riccati equation in \eqref{eqn:value_fun_ss} but using a fixed policy. The policy update step corresponds to \eqref{eqn:policy_update} computing a greedy one-step Nash equilibrium policy with respect to fixed value functions, which relates to the gain expression \eqref{eqn:policy_gain_ss} but using fixed cost matrices.

To the best of our knowledge, the proposed policy iteration algorithm has not been studied and extensively compared with policy gradient methods for general-sum LQDGs.
In the single-player case, the standard policy iteration algorithm coincides with the Gauss-Newton policy gradient method under suitable choice of step size \cite{todorov2005generalized,liao1991convergence}. However, there is a key difference in how the policies are updated in the N-player game setting, which we will explain in detail in the next section. 


\begin{algorithm}
\caption{Policy iteration for N-player general-sum LQDGs with infinite-horizon}\label{alg:policy_iteration_algo}
\begin{algorithmic}[1]
\REQUIRE Stabilizing policies $K_0^i$, system matrices $A$ and $B^i$, cost parameters $Q^i$ and $R^i$, and convergence threshold $\epsilon > 0$.
\STATE \textbf{Initialize:} $K_0^i = (0,\ldots,0)$, $K_1^i = (\infty,\ldots,\infty)$, $i = 1,\ldots,N$, $k = 0$.
\WHILE{$\sum_{i=1}^N \norm{K_{k+1}^{i} - K_k^i} > \epsilon$}
\STATE \textbf{Policy Evaluation:} Compute the value functions of the current policy set for $i = 1,\ldots, N$ by solving the Lyapunov equations
{\small
\begin{align} \label{policyevaluation}
\begin{split}
    P_k^{i} = Q^i + & \rlpar{K_k^{i}}^\top R^i K_k^{i} + \rlpar{A + \sum_{j=1}^{N} B^j K_k^{j}}^\top \\
     & P_k^{i} \rlpar{A + \sum_{j=1}^{N} B^j K_k^{j}}.
\end{split}
\end{align}}
\STATE \textbf{Policy Update:} Update the policy set for $i = 1,\ldots, N$ by solving the coupled equations
{\small
\begin{align}\label{eqn:policy_update}
    \begin{split}
        K_{k+1}^{i} = & -\rlpar{R^i + \rlpar{B^i}^\top P_k^{i} B^i}^{-1} \rlpar{B^i}^\top \\
        & P_k^{i} \rlpar{A  + \sum_{j = 1, j \neq i}^{N} B^j K_{k+1}^{j}},
    \end{split}
\end{align}}
which can be jointly computed by solving \eqref{eqn:ricatti_eqn}.
\STATE $k \leftarrow k+1$
\ENDWHILE
\ENSURE Nash Equilibria $K^{i*}$ $(i = 1,\ldots,N)$
\end{algorithmic}
\end{algorithm}

\begin{figure*}[tb]
{\small
\begin{equation}\label{eqn:ricatti_eqn}
\begin{bmatrix}
R^1+\rlpar{B^1}^\top P_k^1 B^1 & \rlpar{B^1}^\top P_k^1 B^2 & \rlpar{B^1}^\top P_k^1 B^3 & \cdots & \rlpar{B^1}^\top P_k^1 B^N \\
\rlpar{B^2}^\top P_k^2 B^1 & R^2+\rlpar{B^2}^\top P_k^2 B^2 & \rlpar{B^2}^\top P_k^2 B^3 & \cdots & \rlpar{B^2}^\top P_k^2 B^N \\
\rlpar{B^3}^\top P_k^3 B^1 & \rlpar{B^3}^\top P_k^3 B^2 & R^3+\rlpar{B^3}^\top P_k^3 B^3 & \cdots & \rlpar{B^3}^\top P_k^3 B^N \\
\vdots & \vdots & \cdots & \ddots & \vdots \\
\rlpar{B^N}^\top P_k^N B^1 & \rlpar{B^N}^\top P_k^N B^2 & \rlpar{B^N}^\top P_k^N B^3 & \cdots & R^N+\rlpar{B^N}^\top P_k^N B^N
\end{bmatrix}
\begin{bmatrix}
    K_{k+1}^1 \\
    K_{k+1}^2 \\
    K_{k+1}^3 \\
    \vdots \\
    K_{k+1}^N
\end{bmatrix} = 
\begin{bmatrix}
    \rlpar{B^1}^\top P_k^1 \\
    \rlpar{B^2}^\top P_k^2 \\
    \rlpar{B^3}^\top P_k^3 \\
    \vdots \\
    \rlpar{B^N}^\top P_k^N \\
\end{bmatrix} A.
\end{equation}}
\end{figure*}

\section{A Comparison of the Policy Iteration and Gauss-Newton Policy Gradient Algorithms}\label{sec:policy_gradient}
In this section, we first introduce the N-player vanilla/standard policy gradient and natural policy gradient methods, with a brief analysis of the results from \cite{mazumdar2020policy} and \cite{hambly2023policy} in Section \ref{subsec:std_policy_gradient}. Then we extend the Gauss-Newton policy gradient method from the single-player case to the N-player game setting. We show that, unlike the single-player case, the proposed policy iteration algorithm is distinct from the Gauss-Newton policy gradient method in the N-player game setting in Section \ref{subsec:gauss_newton_policy_gradient}.

\subsection{The Vanilla/Standard Policy Gradient and Natural Policy Gradient Methods}\label{subsec:std_policy_gradient}
An extension of the vanilla/standard policy gradient method in the single-player case \cite{fazel2018global} for the N-player game setting has the following form
\begin{align}
    K^i_{k+1} = K^i_k - \eta^i \nabla_{K_k^i} J^i\rlpar{\mathbf{K}_k},
\end{align}
where $\mathbf{K}_k = (K^{1}_k, \ldots, K^{N}_k)$ is a collection of gain matrices of all players with initial gains $\mathbf{K}_0$ and $\eta^i$ is the step size. The work in \cite{mazumdar2020policy} suggests that the vanilla/standard policy gradient method has no guarantees of even local convergence in general-sum infinite-horizon deterministic LQDGs. In contrast, \cite{hambly2023policy} proved the global convergence of the natural policy gradient method to the Nash equilibrium with finite-horizon and stochastic dynamics. The natural policy gradient method presented in \cite{hambly2023policy} is
\begin{equation}
    K^i_{t,k+1} = K^i_{t,k} - \eta \nabla_{K^i_{t}} J^i\rlpar{\mathbf{K}_{k}} \rlpar{\Sigma_{t, \mathbf{K}_{k}}}^{-1},
\end{equation}
which is an N-player extension of the natural policy gradient method introduced in \cite{fazel2018global} and $\Sigma_{t, \mathbf{K}} = \mathbf E [x_{t, \mathbf{K}} x_{t, \mathbf{K}}^\top]$ is the state covariance matrix. However, the natural policy gradient method may fail to converge to a Nash equilibrium for deterministic LQDGs without careful selection of the initial policies and step sizes. To our best knowledge, the general convergence properties of policy gradient algorithms for general-sum LQDGs are not fully understood.
\subsection{The Gauss-Newton Policy Gradient Method}\label{subsec:gauss_newton_policy_gradient}
Our goal in this section is to compare the Gauss-Newton policy gradient method and our algorithm, to highlight and explain the main difference between the two.
For the single-player with stationary feedback policy case, the Gauss-Newton policy gradient method is
\begin{align} \label{eq:GNsinglePlayer}
    K_{k+1} = K_k - \eta \left(R + B^\top P_{k} B\right)^{-1} \nabla J(K_{k}) \Sigma_{k}^{-1},
\end{align}
where the gradient $\nabla J(K_k)$ is equal to
\begin{align*}
    \nabla J(K_k) = 2 \rlpar{ \rlpar{ R+ B^\top P_{k} B } K_k + B^\top P_{k} A } \Sigma_{k}.
\end{align*}
By substituting this in \eqref{eq:GNsinglePlayer} we obtain
\begin{align}\label{eqn:GNPG_single}
    K_{k+1} = (1-2\eta)K_k - 2\eta \left(R + B^\top P_{k} B\right)^{-1} B^\top P_{k} A.
\end{align}

The natural extension of the Gauss-Newton policy gradient method to the N-player game setting would be as
\begin{align*}
K^i_{k+1} = K^i_k - \eta^i \left( R^i + (B^i)^\top P_{k}^{i} B^i \right)^{-1} \nabla_{K_k^i} J^i(\mathbf{K}_k) \Sigma_{k}^{-1}.
\end{align*}
Following the same steps, we obtain
{\small\begin{align} \label{eq:gaussNewtonExpl}
    K_{k+1}^{i} = &(1 - 2\eta^i) K_k^i - 2\eta^i \left( R^i + \rlpar{B^i}^\top P_{k}^{i} B^i \right)^{-1} \nonumber \\
    &\rlpar{B^i}^\top P_{k}^{i} \left( A + \sum_{j=1, j \neq i}^N B^j K_k^j \right).
\end{align}}

In the single-player case, by setting $\eta = \frac{1}{2}$ in \eqref{eqn:GNPG_single}, one can easily verify that the Gauss-Newton policy gradient method is equivalent to the standard policy iteration algorithm \cite{howard:dp}:
\begin{equation}
    K_{k+1} = -\left(R + B^\top P_{k} B\right)^{-1} B^\top P_{k} A.
\end{equation}
For the N-player case, however, this is no longer true. By setting $\eta^i = \frac{1}{2}$ in \eqref{eq:gaussNewtonExpl} we have
\begin{align}
    \begin{split}
        K_{k+1}^{i} ={} &-\left( R^i + \rlpar{B^i}^\top P_{k}^{i} B^i \right)^{-1} \rlpar{B^i}^\top P_{k}^{i} \\
        &\left( A + \sum_{j=1, j \neq i}^N B^j K_k^j \right),
    \end{split}
\end{align}
which is different from the policy update \eqref{eqn:policy_update} of the proposed policy iteration algorithm. In particular, in the Gauss-Newton method, player $i$'s policy update at iteration $k+1$ is defined on the premise that all other players' policy gains $K_k^j$, $j = 1, \ldots, N, j \neq i$ remain fixed at the previous iteration step $k$. In the proposed policy iteration algorithm, however, all players update their policy gains $K_{k+1}^i$, $i = 1, \ldots, N$ simultaneously at iteration $k+1$. As a result, the proposed policy iteration algorithm needs to solve a linear system with $\sum_{i=1}^Nnm_i$ equations and $\sum_{i=1}^Nnm_i$ unknowns as shown in \eqref{eqn:ricatti_eqn} and there has to be a central solver (alternatively, each player can compute their own update independently when the model is assumed to be common knowledge). In the Gauss-Newton policy gradient method, however, each player can compute its own policy gain update.

\section{Numerical Experiments} \label{sec:numerical_experiments}
In this section, we first compare the convergence performance of the proposed policy iteration algorithm, natural policy gradient (Algorithm 1 in \cite{hambly2023policy} extended to the infinite-horizon deterministic case), and Gauss-Newton policy gradient methods under the same experimental setup as in \cite{mazumdar2020policy, hambly2023policy} in \ref{subsec:expA} and \ref{subsec:expB}. Then we compare these algorithms with an additional 1000 random open-loop stable systems with initial policy gain $K^i_0 = (0,\ldots,0)$ that satisfies the conditions in Proposition \ref{prop:LQDG_VI} in \ref{subsec:expC}. 

The model parameters for the numerical results in \ref{subsec:expA} and \ref{subsec:expB} are
\begin{equation*}
    A = 
    \begin{bmatrix}
        0.588 & 0.028 \\ 
        0.570 & 0.056
    \end{bmatrix}, \
    B^1 = 
    \begin{bmatrix}
        1 \\ 
        1
    \end{bmatrix}, \
    B^2 = 
    \begin{bmatrix}
        0 \\ 
        1
    \end{bmatrix},
\end{equation*}

\begin{equation*}
    Q^1 = 
    \begin{bmatrix}
        0.01 & 0 \\
        0 & 1
    \end{bmatrix}, \
    Q^2 = 
    \begin{bmatrix}
        1 & 0 \\
        0 & 0.147
    \end{bmatrix}, \
    R^1 = R^2 = 0.01.
\end{equation*}
We initialize both players' policy $K_0^i$ such that $(A + \sum_{i=1}^N B^i K_0^i)$ is stable. To analyze the convergence performance of all three algorithms under different initial policies, the initial policy gain $K_0^i$ also satisfies $\|K_0^i - K^{i*}\|_2 \leq r$, where $K^{i*}$ denotes a Nash equilibrium of the system, and $r$ is the radius of the ball centered at $K^{i*}$ in which we initialize the policies. A Nash equilibrium of the above system is $K^{1*}=(-0.5134, -0.0439)$ and $K^{2*}=(-0.0525,-0.0114)$. This Nash equilibrium is computed through value iteration introduced in \ref{subsec:value_iteration}. Differently from the normalized error definition in \cite{hambly2023policy}, we define the normalized error of a given pair of policies $(K^1, K^2)$: for $i=1,2$ as
\begin{equation*}
    e_{norm,k} = \frac{\|K_k^1 - K^{1*}\|_2}{\|K^{1*}\|_2} + \frac{\|K_k^2 - K^{2*}\|_2}{\|K^{2*}\|_2},
\end{equation*}
where $k$ is the number of iterations.
\subsection{Faster Convergence Rate of Policy Iteration}\label{subsec:expA}
\begin{figure}[tb] 
    \centering
    \begin{subfigure}[b]{\linewidth}
        \centering
        \includegraphics[width=\linewidth]{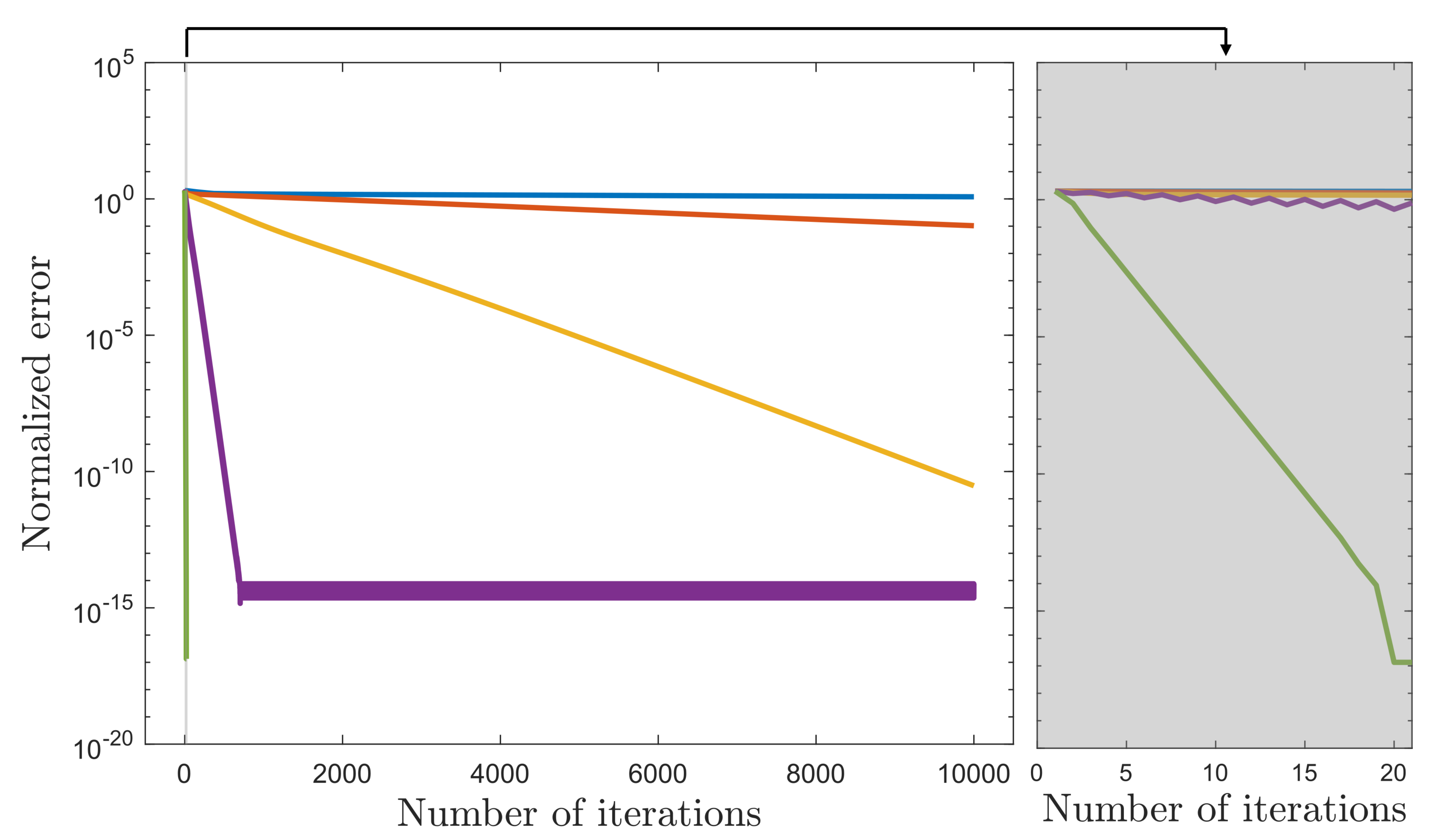}
        \caption{Normalized error vs. number of iterations.}
        \label{fig:num_iterations}
    \end{subfigure}
    \vfill
    \begin{subfigure}[b]{\linewidth}
        \centering
        \includegraphics[width=\linewidth]{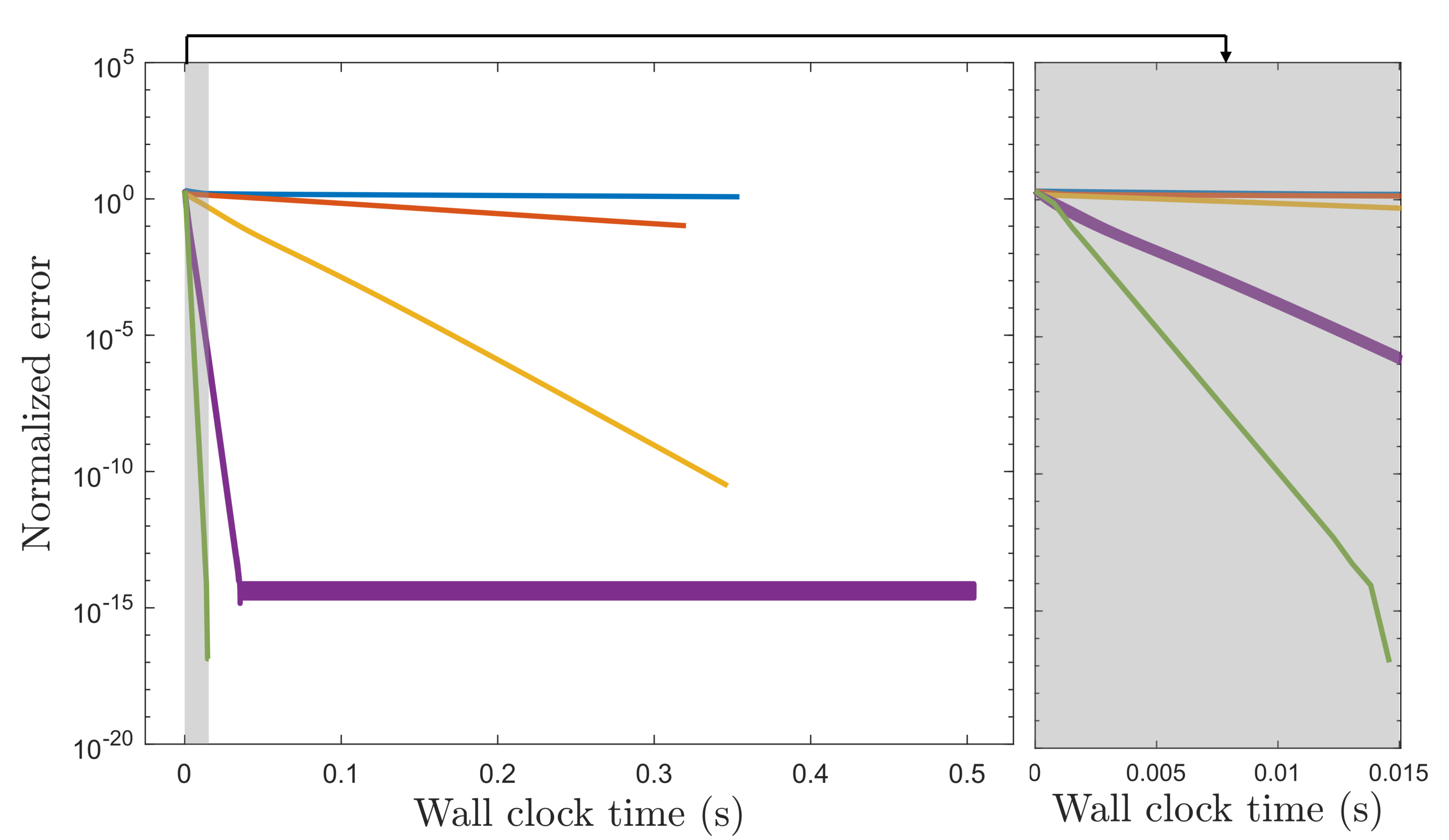}
        \caption{Normalized error vs. elapsed wall clock time.}
        \label{fig:clock_time}
    \end{subfigure}
    \caption{Convergence speed of the proposed policy iteration algorithm (green), Gauss-Newton policy gradient (purple, $\eta^i=0.5$), and natural policy gradient (blue, $\eta^i=10^{-3}$; red, $\eta^i=10^{-2}$; yellow, $\eta^i=10^{-1}$) ($r=0.1$).}
    \label{fig:iteration_time}
\end{figure}

Figure \ref{fig:iteration_time} reveals that all three methods—the proposed policy iteration algorithm, natural policy gradient ($\eta^i=0.1$), and Gauss-Newton policy gradient ($\eta^i=0.5$)—successfully converge to $K^{1*}$ and $K^{2*}$ with the same initial policy gains $K_0^1 = (-0.4266, -0.0938)$ and $K_0^2 = (0.0342, -0.0612)$ ($r=0.1$). It is evident that the proposed policy iteration algorithm converges to the Nash equilibrium with much fewer iterations and shorter computational time compared with the other two policy gradient methods in this specific example. The natural policy gradient method fails to converge to the Nash equilibrium if the step size is not carefully selected. Although the Gauss-Newton policy gradient method converges to the Nash equilibrium, it requires significantly more iterations and computational time than the proposed policy iteration algorithm. However, comparing the general theoretical convergence rates of the proposed policy iteration algorithm and policy gradient methods remains open.
\subsection{Convergence Performance of Policy Iteration from a Distant Initial Policy}\label{subsec:expB}
\begin{figure}[htbp] 
    \centering
    \begin{subfigure}[b]{\linewidth}
        \centering
        \includegraphics[width=0.95\linewidth]{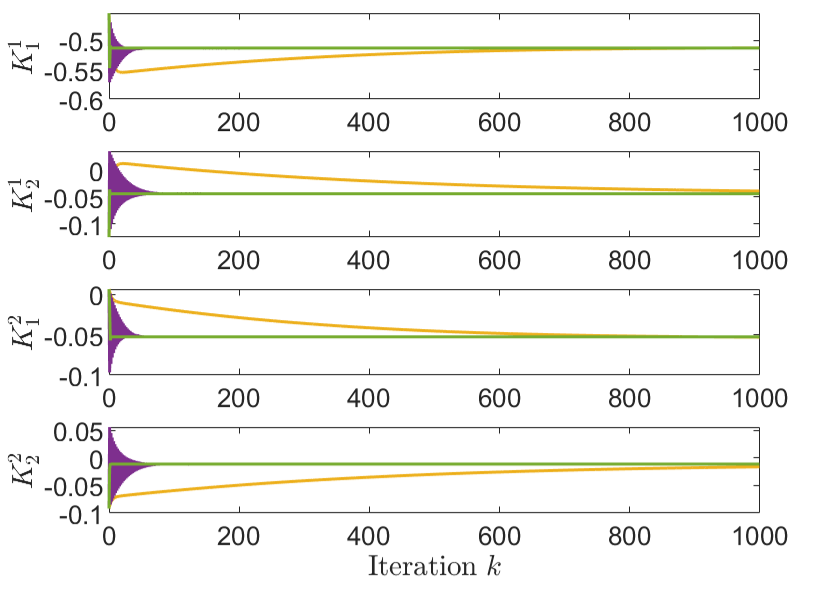}
        \caption{An instance ($r = 0.1$) that all three algorithms converge to $K^{1*}$ and $K^{2*}$.}
        \label{fig:K_traj_small_radius}
    \end{subfigure}
    \vfill
    \begin{subfigure}[b]{\linewidth}
        \centering
        \includegraphics[width=0.95\linewidth]{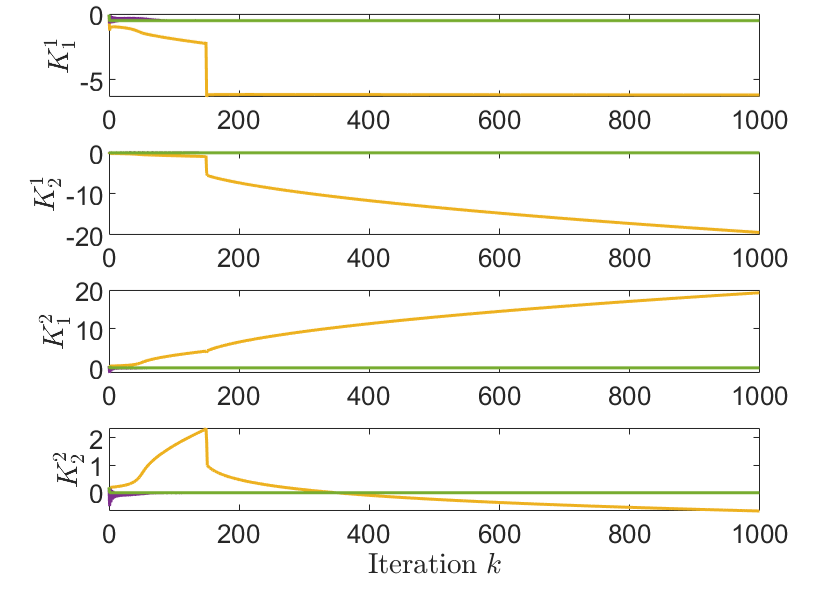}
        \caption{An instance ($r=0.5$) that the natural policy gradient algorithm fails to converge to $K^{1*}$ and $K^{2*}$.}
        \label{fig:K_traj_big_radius}
    \end{subfigure}
    \caption{Convergence performance of the proposed policy iteration algorithm (green), Gauss-Newton policy gradient (purple, $\eta^i=0.5$), and natural policy gradient (yellow, $\eta^i=10^{-1}$) methods under different initial policy gains.}
    \label{fig:conv_perfm_init_policy}
\end{figure}

Figure \ref{fig:conv_perfm_init_policy} illustrates the convergence performance of the proposed policy iteration algorithm, natural policy gradient, and Gauss-Newton policy gradient methods as the initial policy gains transition from a smaller to a larger neighborhood around $K^{1*}$ and $K^{2*}$. Figure \ref{fig:K_traj_small_radius} presents an instance where all three methods with the same initial policy gains $K_0^1$ and $K_0^2$ as in \ref{subsec:expA} converge to $K^{1*}$ and $K^{2*}$. Conversely, Figure \ref{fig:K_traj_big_radius} shows a different scenario where, under a pair of more distant initial policy gains $K_0^1 = (-0.0543, 0.1541)$ and $K_0^2 = (0.4066, 0.1867)$ ($r=0.5$) from $K^{1*}$ and $K^{2*}$, the natural policy gradient method does not converge to $K^{1*}$ and $K^{2*}$ as the other two methods. It is evident that the convergence performance of the proposed policy iteration algorithm is less sensitive to changes in initial policy compared to policy gradient methods. It is possible that the convergence performance of the policy gradient methods can be improved by tuning the step size. However, this process for each case can be laborious and time-consuming.
\subsection{Convergence Performance of Policy Iteration for Additional Problem Instances}\label{subsec:expC}
We now present results for additional problem instances by randomly generating system parameters ($A, B^i, Q^i, \text{and} \ R^i$). The entries of these parameters were independently drawn from a standard normal distribution. The dynamics matrix is scaled to be open-loop stable, and all cost parameters are made positive definite. In each instance, value iteration was used to compute a Nash equilibrium (if it converges) for that game setting. Then the policy iteration and policy gradient methods with the same (stabilizing) initial policy gain $K_0^i = (0,\ldots,0)$ were used to compute a solution for the same game setting. We examine whether and how fast the methods converge to the Nash equilibrium computed using value iteration.

Table \ref{tab:random_sys_table} shows the convergence performance of the proposed policy iteration algorithm, natural policy gradient, and Gauss-Newton policy gradient methods for 1000 randomly-generated problem instances for both two and four players. In all experiments, the convergence of the proposed policy iteration algorithm is significantly faster and more reliable than the policy gradient methods. Moreover, the convergence performance of the proposed policy iteration algorithm is less sensitive to a change in the number of players from two to four. 
These results provide additional evidence that the proposed policy iteration algorithm outperforms policy gradient methods, especially since there is no need to tune the step size for each instance.
\begin{table}[htbp]
    \setlength\tabcolsep{0pt}
    \centering
    \caption{Convergence performance of the policy iteration (PI) algorithm, natural policy gradient (NPG, $\eta^i = 0.1$), and gauss-Newton policy gradient (GNPG, $\eta^i = 0.5$) methods for $1000$ random systems.}
    \begin{tabular}{|c|c|c|c|c|}
        \hline
        & \multicolumn{2}{|c|}{$n = 4, m^i = 2, N = 2$} & \multicolumn{2}{|c|}{$n = 4, m^i = 2, N = 4$} \\
        \hline
        & \ convergent \ & \ average number \ & \ convergent \ & \ average number \ \\
        & cases & of iterations & cases & of iterations \\
        \hline
         \ NPG \ & $5$ & $139$ & $0$ & N/A \\
         \ GNPG \ & $880$ & $161$ & $8$ & $1259$ \\
         \ PI \ & $960$ & $7$ & $945$ & $7$ \\
         \hline
    \end{tabular}
    \label{tab:random_sys_table}
\end{table}

In all our empirical studies, the proposed policy iteration algorithm converges at a much higher speed than the policy gradient methods. Similar numerical results have also been shown in \cite{nortmannNash2024}. This may be due to the fact that the proposed policy iteration algorithm takes into account the policies of the other players at iteration $k+1$. The update in the policy gradient methods considers only the policy of the other players at iteration $k$. Thanks to this additional information, the policy iteration algorithm adjusts the update of the policy and therefore, we believe, avoids possible overshooting.

\section{Conclusions} \label{sec:conclusions}
We proposed a policy iteration algorithm for the infinite-horizon N-player general-sum deterministic LQDGs and compare it to policy gradient methods. We demonstrated that the proposed policy iteration algorithm is distinct from the Gauss-Newton policy gradient method in the N-player game setting, in contrast to the single-player case where they are equivalent under suitable choice of step size. 
In all our numerical experiments, the proposed policy iteration algorithm converges to the same Nash equilibrium as value iteration with fewer iterations and less computational time compared to policy gradient methods. Furthermore, the convergence performance of the proposed policy iteration algorithm is less sensitive to the initial policy and changes in the number of players. 


\bibliographystyle{IEEEtran}
\bibliography{bibliography.bib}

\end{document}